\documentclass[a4paper,11pt]{amsart}

\usepackage{latexsym}
\usepackage{amsmath,amsfonts,amssymb}
\usepackage[english]{babel}
\usepackage[all]{xy}
\usepackage{graphicx}

\newcommand{\al}{\alpha}
\newcommand{\be}{\beta}
\newcommand{\ga}{\gamma}
\newcommand{\ep}{\epsilon}
\newcommand{\de}{\delta}
\newcommand{\si}{\sigma}

\newcommand{\la}{\lambda}

\newcommand{\om}{\omega}
\newcommand{\De}{\Delta}
\newcommand{\La}{\Lambda}
\newcommand{\Ga}{\Gamma}
\newcommand{\Om}{\Omega}

\newcommand{\N}{\mathbb{N}}
\newcommand{\Z}{\mathbb{Z}}
\newcommand{\Q}{\mathbb{Q}}
\newcommand{\R}{\mathbb{R}}
\newcommand{\B}{\mathbb{B}}

\newcommand{\M}{\mathcal{M}}
\newcommand{\A}{\mathcal{A}}
\newcommand{\G}{\mathcal{G}}
\newcommand{\W}{\mathcal{W}}
\newcommand{\Pj}{\mathbb{P}}
\newcommand{\Rees}{\mathfrak{R}}
\newcommand{\field}{\Bbbk}
\newcommand{\lig}{\mathfrak{g}}
\newcommand{\lih}{\mathfrak{h}}
\newcommand{\lit}{\mathfrak{t}}

\newcommand{\bn}{{\bf n}}
\newcommand{\perm}{\mathfrak{S}}

\newcommand{\coal}{\al^\vee}

\newcommand{\tal}{{\tilde{\al}}}
\newcommand{\oal}{\overline{\al}}
\newcommand{\osi}{\overline{\si}}
\newcommand{\ola}{\overline{\la}}
\newcommand{\omu}{\overline{\mu}}
\newcommand{\ox}{\overline{x}}
\newcommand{\oz}{\overline{z}}
\newcommand{\tom}{\tilde{\om}}
\newcommand{\tPhi}{\tilde{\Phi}}
\newcommand{\tDe}{\tilde{\De}}

\newcommand{\tW}{\tilde{W}}
\newcommand{\tV}{\tilde{V}}
\newcommand{\tA}{\tilde{A}}
\newcommand{\tT}{\tilde{T}}
\newcommand{\tS}{\tilde{S}}
\newcommand{\tbe}{\tilde{\be}}
\newcommand{\cG}{{\tilde{G}}}
\newcommand{\cH}{{\tilde{H}}}
\newcommand{\cP}{{\tilde{P}}}

\newcommand{\tlb}{\mathcal{O}(1)}
\newcommand{\lb}{\mathcal{L}}

\newcommand{\ch}{\mathop{\rm ch}\nolimits}
\newcommand{\Id}{\mathop{\rm Id}\nolimits}
\newcommand{\Pic}{\mathop{\rm Pic}\nolimits}
\newcommand{\End}{\mathop{\rm End}\nolimits}
\newcommand{\Hom}{\mathop{\rm Hom}\nolimits}

\newcommand{\he}{\mathop{\rm ht}\nolimits}
\newcommand{\img}{\mathop{\rm Im}\nolimits}

\newcommand{\lex}{\mathop \leq_{\rm lex}\nolimits}

\newcommand{\typeBC}{{\textrm{\textsf{BC}}}}
\newcommand{\typeB}{{\textrm{\textsf{B}}}}
\newcommand{\typeF}{{\textrm{\textsf{F}}}}

\newtheorem{definition}{Definition}[section]
\newtheorem{theorem}[definition]{Theorem}
\newtheorem{corollary}[definition]{Corollary}

\newtheorem{proposition}[definition]{Proposition}
\newtheorem{lemma}[definition]{Lemma}
\newtheorem{example1}[definition]{Example}

\newtheorem{remark1}[definition]{Remark}

\title{The ring of sections of a complete symmetric variety}

\author{Rocco Chiriv\`\i\ and Andrea Maffei}

\thanks{The first named author was partially supported by a grant of the Department of Mathematics of the University of Rome ``La Sapienza''}

\keywords{Complete symmetric variety, Picard group, PRV conjecture, standard monomial theory, rational singularities}

\begin{document}

\begin{abstract}
We study the ring of sections $A(X)$ of a complete symmetric variety $X$, that is of the wonderful completion of $G/H$ where $G$ is an adjoint semi-simple group and $H$ is the fixed subgroup for an involutorial automorphism of $G$. We find generators for $\Pic(X)$, we generalize the PRV conjecture to complete symmetric varieties and construct a standard monomial theory for $A(X)$ that is compatible with $G$ orbit closures in $X$. This gives a degeneration result and the rational singularityness for $A(X)$.
\end{abstract}

\maketitle

\section*{Introduction}

The aim of this paper is to explicitly adapt the Littelmann standard monomial theory for 
flag varieties in \cite{L1},\cite{L2} and \cite{LLM}, to complete symmetric varieties as 
constructed by De Concini and Procesi in \cite{CP} in characteristic zero and by 
De Concini and Springer for arbitrary characteristic in \cite{CS}.

We review briefly such completions. Let $G$ be an adjoint semi-simple group, 
let $H$ be the fixed subgroup for an involutive automorphism $\si$ of $G$ and consider 
the affine variety $G/H$, called \emph{symmetric variety}. 
De Concini and Procesi in \cite{CP} show that there exists 
an irreducible representation $V$ and a 
line $r$ of $V$ such that the stabilizer of $r$ is $H$. 
They define the complete symmetric variety $X$ as the closure of the orbit $G\cdot r$
in $\Pj (V)$.

This compactification of $G/H$ is a \emph{wonderful} $G$ variety in the 
sense of Luna (see Proposition \ref{pcompactification}): $X$ is smooth; $X\setminus G\cdot r$ is a divisor with normal crossing and smooth irreducible components $S_1,\ldots,S_\ell$; the closures of the $G$ orbits in $X$ are $X_I=\cap_{i\in I}S_i$ where $I\subset\{1,\ldots,\ell\}$.

Many properties of these varieties have been studied.
The ones interesting for us here are the description of the Picard group
of $X$ as a sublattice of the lattice of weights given in 
\cite{CP} and \cite{CS} 
and the description of $H^0(X,\lb)$ as a $G$-module for all $\lb \in \Pic(X)$. In particular in \cite{CP} it is proved that $H^0(X,\lb_{\al_i-\si(\al_i)})\simeq\field s_i$ where $\al_i$ is a simple root for a suitable basis of the root system of $G$ and $s_i$ is a $G$ invariant section whose divisor is $S_i$.

In this paper we construct a standard monomial theory for the ring 
$$
A(X) = \bigoplus _{\lb \in \Pic(X)} H^0(X,\lb).
$$
We call this ring the \emph{ring of sections} of $X$. One should think to this ring as a variant of the multicone over a flag variety.

We give first an explicit description of a basis of $\Pic(X)$ (see Theorem 
\ref{tspherical}). As a first application of this result we prove a generalization 
of the Parthasarathy-Ranga Rao-Varadarajan conjecture (PRV) to 
complete symmetric varieties (Theorem \ref{tsprv}). The motivation for this generalization
is that it should give the combinatorial side of the surjectivity of the 
multiplication map $H^0(X,\lb)\otimes H^0(X,\lb')\rightarrow H^0(X,\lb\otimes\lb').$
In the so-called group case, i.e. the involution 
$\si:G\times G\rightarrow G\times G$ with $\si(g_1,g_2)=(g_2,g_1)$, 
the surjectivity has been proved 
by Kannan in \cite{K} but in general this remains an open problem.

Then we pass to the construction of a standard monomial theory for the ring 
$A(X)$. Let $P$ be the parabolic subgroup of $G$ such that $G/P$ is the unique closed orbit in $X$. Given an LS path $\pi$ of shape $\la\in\Pic^+(G/P)$, Littelmann defines a section $p_\pi$ in $H^0(G/P,\lb_\la)$. 
We lift $p_\pi$ to a section $x_\pi$ over $X$ taking into account the description 
of $H^0(X,\lb_\la)$. The building blocks of our monomials are given by the sections 
$x_\pi$, where $\pi$ runs over LS paths of shape $\theta$ for $\theta$ generator of 
$\Pic(X)$, and by the sections $s_1,\ldots,s_\ell$. Then we define a notion 
of standardness for these monomials and introduce a variant of the lexicographic order. 
Our standard monomial theory strictly mimes that of $G/P$, indeed the relations are the 
same up to ``bigger'' terms that vanish ``more'' on the divisors $S_1,\ldots,S_\ell$ 
(Theorem \ref{tsmt}). Further this standard monomial theory is compatible with the $G$ orbit closures (Corollary \ref{cgvarietybasis}).

As a consequence of this standard monomial theory we construct a flat deformation that degenerates $A(X)$ to $\tA(G/P)\otimes\field[s_1,\ldots,s_\ell]$, where $\tA(G/P)$ is the coordinate ring of a multicone over $G/P$ corresponding to the sublattice $\Pic(X)$ 
of $\Pic(G/P)$.
So we use this to prove that $A(X)$ has rational singularities. 
Moreover if we fix a line bundle $\lb$ then also the ring 
$A_\lb=\oplus_n H^0(X,\lb^{\otimes n})$ has rational singularities. This is well known and we have included it here since it was impossible for us to find it in the literature. 

Although in all this paper we take $\field$ to have characteristic zero, all the results are valid in every characteristic except the rational singularityness proofs that use the flat degeneration and the quotient by a reductive group to pass from $A(G/P)$ to $A(X)$.

The paper is organized as follows. In Section \ref{spreliminary} we recall 
all preliminaries results about complete symmetric varieties to be used in the 
sequel. In Section \ref{spicard} we see the description of the generators 
of the Picard group of a symmetric variety giving the application to the 
PRV for symmetric varieties. Finally in Section \ref{ssmt}, after a short 
review of Littelmann LS paths and related results, we construct our 
standard monomial theory. We finish the section proving the degeneration result and showing that the ring $A(X)$ and the cone ring $A_\lb(X)$ have rational singularities.

The authors would like to thank C. De Concini who suggested us to work 
on a standard monomial theory for the ring $A(X)$. We want also to thank him for many 
useful conversations on this problem.

\section{Preliminary Results on Complete Symmetric Varieties}\label{spreliminary}

In this section we collect all preliminary results for the sequel setting up
notation and reviewing the construction of the wonderful
compactification of $G/H$ (for details see \cite{CP} and \cite{CS}).

Let $G$ be an adjoint semi-simple group defined over an algebraically closed field $\field$ of
characteristic zero, and let $\si$ be an involutorial automorphism of $G$.
Denote by $H$ the subgroup of fixed points of $\si$ in $G$. The involution
$\si$ induces a linear map, still denoted by $\si$, on the Lie algebra $\lig$
of the group $G$. We denote by
$\lih$ the Lie algebra of the reductive group $H$; notice that $\lih$ is
exactly the $+1$ eigenspace of $\si$ on $\lig$. If $T$ is a $\si$ stable
torus of $G$ and $\lit$ its Lie algebra, we decompose $\lit$ as
$\lit_0\oplus\lit_1$ with $\lit_0$ the $+1$ eigenspace of $\si$ and $\lit_1$
the $-1$ eigenspace. Notice that $\lit_0$ is the Lie algebra of $T^\si$ while
$\lit_1$ is the Lie algebra of the torus $T_1=\{t\in T\ |\ \si(t)=t^{-1}\}$;
we call this latter torus \emph{anisotropic}. Recall that any $\si$ stable
torus is contained in a maximal torus of $G$ which is itself $\si$ stable.
We fix such a $\si$ stable maximal torus $T$ for which $\dim T_1$ is maximal
and denote this dimension by $\ell$, calling it the \emph{rank} of the
symmetric variety $G/H$.

Now let $\Phi\subset\lit^*$ be the root system of $\lig$ and denote still by
$\si$ the induced map on $\lit^*$. Observe that $\si$ preserves the killing
form on $\lit$ and on $\lit^*$. Let $\Phi_0=\{\al\in\Phi\ |\ \si(\al)=\al\}$
and $\Phi_1=\Phi\setminus\Phi_0$.  We can choose the set of positive roots
$\Phi^+$ in such a way that $\si(\al)\in\Phi^-$ for all root
$\al\in\Phi^+\cap\Phi_1$. Let $\De$ be the basis defined by $\Phi^+$ and put
$\De_0=\De\cap\Phi_0$, $\De_1=\De\cap\Phi_1$. The action of the involution
$\si$ on the set of roots admits the following descriptions. There exists an
involutive bijection $\osi:\De_1\rightarrow\De_1$ such that for every
$\al\in\De_1$ we have
$$
\si(\al)=-\osi(\al)-\be_\al
$$
where $\be_\al$ is a non negative linear combination of roots in $\De_0$,
moreover $\be_{\osi(\al)}=\be_\al$.
Further $\si(\al)=-w_{\De_0}\osi(\al)$ if $\al\in\De_1$, where $w_{\De_0}$ is the longest element of the Weyl group of the root system with base $\De_0$.

We introduce here a particular behavior
of a simple root: we say that $\al\in\De_1$ is an \emph{exceptional} root if
$\osi(\al)\neq\al$ and $(\al,\si(\al))\neq0$, where $(\cdot,\cdot)$ is the
Killing form. Notice that $\osi(\al)$ is exceptional if $\al$ is. Moreover
the compactification $X$ we are going to construct below is said to be exceptional if
there exist exceptional roots.

Denote by $\La\subset\lit^*$ the set of integral weights of $\Phi$;
clearly $\si$ acts also on this set. Let $\La^+$ be the set of dominant
weights with respect to $\De$ and let
$\om_\al$ be the fundamental weight dual to the simple coroot $\coal$ for
$\al\in\De$. A simple computation, using the $\si$ invariance of the Killing
form, shows that
$$
\si(\om_\al)=-\om_{\osi(\al)}
$$
for every $\al\in\De_1$. Moreover any integral weight $\la$ such that
$\si(\la)=-\la$ is of the form
$$
\la=\sum_{\al\in\De_1}n_\al\om_\al
$$
with integer coefficients such that $n_\al=n_{\osi(\al)}$; we denote by
$\La_1$ the set of such weights that we call \emph{special}. Further if
$n_\al\neq0$ for every root $\al\in\De_1$ then we say that the special
weight $\la$ is \emph{regular}.

Other notation we will use. We index $\De$ as $\al_1,\ldots,\al_n$ and we denote by $\om_i$ the fundamental weight $\om_{\al_i}$. Moreover we define a map $\si$ on the set of indeces $\{1,\ldots,n\}$ in such a way that $\osi(\al_i)=\al_{\si(i)}$ if $\al_i\in\De_1$ and $\si(i)=i$ if $\al_i\in\De_0$. So we can write the action of $\si$ on the fundamental weights of $\De_1$ simply as $\si(\om_i)=-\om_{\si(i)}$.

Now we come to the basic construction of the compactification of $G/H$. Consider a simply connected covering $\pi:\cG\rightarrow G$ and the induced involutorial automorphism $\si:\cG\rightarrow\cG$. For a subgroup $A$ of $G$ we denote by $\tilde{A}$ the subgroup $\pi^{-1}(A)$ of $\cG$. Notice that $\cH=\pi^{-1}(H)$ contains $(\cH)^0=(\cG)^\si$, the fixed point group in $\cG$, as the identity component. If $V$ is a $\cG$ module, we define the $\si$ twisted $\cG$ module $V^\si$ as the vector space $V$ with action $g\cdot v=\si(g)v$. Notice that if $\la$ is a dominant special weight then the dual of the irreducible $\cG$ module of highest weight $\la$ is isomorphic to $V^\si_\la$. Let $h_\la$ be such isomorphism considered as an element of $V_\la\otimes V_\la$. Notice that $G$ acts on the projective space over any $\cG$ module.

Assume now that $\la$ is a regular special dominant weight. In \cite{CP} it is proved that: (i) the stabilizer of the line $\field h_\la$ in $\cG$ is $\cH$ and (ii) the stabilizer of the line $\field h_\la$ in $G$ is $H$.

Consider the $\cG$ decomposition $V_{2\la}\oplus V'$ of the tensor product $V_\la\otimes V_\la$, where $V'$ is sum of highest weight modules $V_\mu$ with $\mu<\la$ in the dominant order, and let $p:V_{2\la}\oplus V'\rightarrow V_{2\la}$ be the $\cG$ equivariant projection. Notice that $p(h_\la)$ is non zero and let $r_\la$ be its class in $\Pj(V_{2\la})$.
Now we define the \emph{compactification} $X$ of $G/H$ as the closure in $\Pj(V_{2\la})$ of the orbit $G\cdot r_\la$. 
Let $P$ be the parabolic subgroup of $G$ stabilizing the line $\field v_{\lambda} \in \Pj(V_{\lambda})$ spanned by a heighest weight vector $v_\la$.
The following proposition from \cite{CP} describes the structure of the compactification.
\begin{proposition}\label{pcompactification}\emph{(Theorem 3.1 in \cite{CP})}\\
1) $X$ is a smooth projective $G$ variety;\\
2)$X\setminus G\cdot r_\la$ is a divisor with normal crossing and smooth irreducible components $S_1,\ldots,S_\ell$;\\
3) the $G$ orbits of $X$ correspond to the subsets of the indexes $1,2,\ldots,\ell$ so that the orbit closures are the intersections $S_{i_1}\cap S_{i_2}\cap\cdots\cap S_{i_k}$, with $1\leq i_1,\ldots,i_k\leq\ell$.\\
So $X$ is a wonderful $G$ variety in the sense of Luna. Moreover\\
4) the unique closed orbit $Y\doteq\cap_{i=1}^\ell S_i$ is isomorphic to the flag variety $G/P$;\\
5) $X$ is independent on the choice of the regular special weight $\la$ up to $G$ equivariant isomorphism.
\end{proposition}

We go on constructing some line bundles on the variety $X$. Let $\la$ be a dominant weight of $\lig$ such that $\Pj(V_\la)$ contains a line $r$ invariant for $\cH$. Consider the map
$$
G/H\ni gH\mapsto g\cdot r\in\Pj(V_\la).
$$
One can show that this induces a projection
$$
\psi_\la:X\rightarrow\Pj(V_\la).
$$
Now let $\tlb$ be the tautological line bundle on $\Pj(V_\la)$ and define the line bundle $\lb_\la$ on $X$ as $\psi_\la^*\tlb$. If we restrict $\lb_\la$ on $G/P\simeq\cG/\cP\simeq Y\hookrightarrow X$ we have the usual line bundle $\cG\times_\cP\field_{-\la}$ corresponding to $\la$ in the identification of $\Pic(\cG/\cP)$ with a sublattice of the weight lattice $\La$. Moreover we have
\begin{proposition}\label{pinjectivepic}\emph{(Proposition 8.1 in \cite{CP})}
The map $\Pic(X)\rightarrow\Pic(Y)$ induced by the inclusion is injective.
\end{proposition}
So we can identify $\Pic(X)$ with a sublattice of the weight lattice. Further the line bundles constructed above account for all line bundles since we have
\begin{proposition}\label{tpic}\emph{(Lemma 4.6 in \cite{CS})}
$\Pic(X)$ corresponds to the lattice generated by the dominant weights $\la$ such that $\Pj(V_\la)^\cH$ is non void.
\end{proposition}
We come to the analysis of such weights. Call a dominant weight $\la$ \emph{spherical} if $V_\la^{\cH^0}\neq0$. It is easy to see that a spherical weight must be special. On the contrary the double of any dominant special weight is spherical. Moreover if $\cH^0$ has no non trivial characters then $\Pic(X)\cap\La^+$ is exactly the set of spherical weights. In general we have
\begin{proposition}\label{ppic}\emph{(Theorem 4.8 in \cite{CS})}
$\Pic(X)$ is generated by the spherical weights and the fundamental weights corresponding to the exceptional roots.
\end{proposition}
We recall a characterization of the spherical weights due to Helgason (see \cite{H}, \cite{S} or \cite{V}). For a root $\al$ let $\oal$ be its restriction to $\lit_1$.
Then
\begin{proposition}\label{thelgason}\emph{(Theorem 3 in \cite{V})}
A weight is spherical if and only if it is special and $(\mu,\oal)/(\oal,\oal)$ is an integer for all root $\al$ such that $(\oal,\oal)\neq0$.
\end{proposition}

Now we introduce a filtration on the spaces of global sections $H^0(X,\lb_\la)$ by the order of vanishing on the $G$ stable divisors $S_1,\ldots,S_\ell$. For a root $\al$ define $\tal\doteq\al-\si(\al)$  and notice that the set $\De_1=\{\al_1,\ldots,\al_\ell,\al_{\ell+1},\ldots,\al_r\}$ can be indexed in such a way that $\tal_i=\al_i-\si(\al_i)$ are different for $i=1,\ldots,\ell$ (and $\tDe=\tDe_1=\{\tal_1,\ldots,\tal_\ell\}$). Then, up to reindexing the $G$ stable divisors, we have
\begin{proposition}\label{pdivisor}\emph{(Corollary 8.2 in \cite{CP})}
There exists a unique up to scalar $G$ invariant section $s_i\in H^0(X,\lb_{\tal_i})$ whose divisor is $S_i$.
\end{proposition}
For a $\ell$-tuple $\bn=(n_1,\cdots,n_\ell)$ of non negative integers, the multiplication by $s^\bn\doteq \Pi_i s_i^{n_i}$ gives a linear map
$$
H^0(X,\lb_{\la-\sum n_i\tal_i})\rightarrow H^0(X,\lb_\la).
$$
Let $F_\la(\bn)$ be the image of this map. We order $\N^\ell$ by setting $(n_1,\ldots,n_\ell)\geq(n'_1,\ldots,n'_\ell)$ if $n_i\geq n'_i$ for $i=1,\ldots,\ell$. Clearly $F_\la(\bn')\subset F_\la(\bn)$ if and only if $\bn'\geq\bn$. We have the following theorem
\begin{proposition}\label{pfiltration}\emph{(Theorem 5.10 in \cite{CP})}
Let $\la\in\Pic(X)$. If $\la-\sum n_i\tal_i$ is dominant, then
$$
F_\la(\bn)/(\sum_{\bn'>\bn}F_\la(\bn'))=H^0(G/P,\lb_{\la-\sum n_i\tal_i}).
$$
Otherwise both sides are $0$. In particular $H^0(X,\lb_\la)\neq0$ if and only if there exists a dominant weight $\mu$ and a $\ell$-tuple of non negative integers $(n_1,\ldots,n_\ell)$ such that $\la=\mu+\sum n_i\tal_i$.
\end{proposition}
As a direct consequence we have
\begin{corollary}\label{csurjective}
Let $\la$ be a dominant weight in $\Pic(X)$. Then the map
$$
H^0(X,\lb_\la)\rightarrow H^0(G/P,\lb_\la)
$$
induced by inclusion, is surjective.
\end{corollary}

We finish this review of preliminary results introducing the restricted root system. The results we state here are proved in \cite{R}. Denote by $\tPhi$ the set $\{\tal\ |\ \al\in\Phi\}$. This is a root system in the space $E_1\doteq\La_1\otimes\R$ with base $\tDe=\{\tal_1,\ldots,\tal_\ell\}$. We call its Weyl group $\tW$ the \emph{restricted Weyl group}. Consider the following subgroups of the Weyl group $W$ of $\lig$, $W_0=\{w\in W\ |\ w(E_1)\subset E_1\}$ and $W_1=\{w\in W\ |\ w_{|E_1}=\Id_{E_1}\}$. Then
\begin{proposition}\label{prestricted}\emph{(Lemma 4.1 in \cite{R})}
The restriction map $W_0\ni w\mapsto w_{|E_1}\in\End_\R(E_1)$ induces an isomorphism of $W_0/W_1$ with the Weyl group $\tW$ of the root system $\tPhi$.
\end{proposition}
Let $\Om_1=\{\mu\in\La_1\ |\ \mu\textrm{ is integral on }\tPhi^\vee\}$ and notice that $\Om_1$ can be identified with the set of integral weights of the root system $(\tPhi,E_1)$.

\section{The Spherical Weights and the PRV Conjecture}\label{spicard}

In this section we complete the description of the spherical weights. Using this description we prove a version of the Parthasaraty-Ranga Rao-Varadarajan conjecture (PRV) for complete symmetric varieties. We begin with two preliminary lemmas.
\begin{lemma}\label{ldomspherical}
If $\la$ is a dominant weight in the lattice generated by the spherical weights then $\la$ is spherical.
\end{lemma}
\begin{proof}
This is clear from the Helgason criterion in Proposition \ref{thelgason}.
\end{proof}
\begin{lemma}\label{lexceptionalval}
Let $\al\in\De_1$ be an exceptional root. Then
$\langle\si(\al),\al^\vee\rangle=1$.
\end{lemma}
\begin{proof}
We have that $\si(\al)$ is not supported in $\al$ since $\si(\al)=-w_{\De_0}\osi(\al)$ and $\osi(\al)\neq\al$. Hence $(\si(\al),\al)\geq0$.

We know that $|\langle\si(\al),\al^\vee\rangle|=|2(\si(\al),\al)/(\al,\al)|\neq2$ since $\si$ preserves the Killing form. Moreover $(\si(\al),\al)\neq0$,
$\si(\al)\neq-\al$ being $\al$ exceptional, so
$\langle\si(\al),\al^\vee\rangle=\pm1$ and by $(\si(\al),\al)\geq0$
we conclude $\langle\si(\al),\al^\vee\rangle=1$.
\end{proof}
The following theorem gives a new description of the spherical weights.
(Recall that we have indexed the set $\De_1$ in such a way that
$\al_i-\si(\al_i)$ are different for $i=1,\ldots,\ell$.)
\begin{theorem}\label{tspherical}
The lattice generated by the spherical weights is the lattice $\Om_1$ of
integral weights of $\tPhi$. Moreover if we set:
$$
\tom_i=\left\{
\begin{array}{ll}
\om_i  & \textrm{if } \si(i)=i,\ \si(\al_i)\neq-\al_i,\\
2\om_i & \textrm{if } \si(\al_i)=-\al_i,\\
\om_i+\om_{\si(i)} & \textrm{if } \si(i)\neq i.\\
\end{array}
\right.
$$
then $\Om_1$ is generated by $\tom_1,\ldots,\tom_\ell$ and $\langle\tom_i,\tal_j^\vee\rangle=c_i\de_{i,j}$ with $c_i=1$ if $2\tal_i\not\in\tPhi$ and $c_i=2$ otherwise. In particular if $\tPhi$ is reduced then $\tom_1,\ldots,\tom_\ell$ are the fundamental weights dual to $\tal_1^\vee,\ldots,\tal_\ell^\vee$.
\end{theorem}
\begin{proof}
In order to describe the spherical weights we use the Helgason criterion
in Proposition \ref{thelgason}.
For each root $\al\in\Phi$ such that $(\oal,\oal)\neq0$ we have
$$
\frac{(\mu,\oal)}{(\oal,\oal)}=\frac{2(\mu,\tal)}{(\tal,\tal)}=
\langle\mu,\tal^\vee\rangle
$$
since $\oal=\frac{1}{2}\tal$. The first claim follows.

Now let $c_\al\doteq2-\langle\si(\al),\al^\vee\rangle$ for $\al\in\Phi$. Observe that $c_\al$ is a non negative integer less or equal to $4$ since $\si$ preserves the Killing form. Moreover $c_\al=0$ if and only if $\si(\al)=\al$, and $c_\al=4$ if and only if $\si(\al)=-\al$. Hence, in particular, $(\oal,\oal)\neq0$ implies $c_\al\neq0$.

Fix a root $\al$ such that $(\oal,\oal)\neq0$ and let $\mu$ be a special dominant weight. We have
$$
\frac{(\mu,\oal)}{(\oal,\oal)}=\frac{2\langle\mu,\al^\vee\rangle}{c_\al}.
$$
In \cite{CP} it is proved that $2\mu$ is a spherical weight (see also the discussion above Proposition \ref{pcompactification}). So $(\mu,\oal)/(\oal,\oal)=a/2$ with $a\in\Z$, hence we have the integral equality
$4\langle\mu,\al^\vee\rangle=a c_\al$.

There are two cases: (i) $\si(\al)\neq-\al$, so $4\nmid c_\al$, hence $2|a$ and $(\mu,\oal)/(\oal,\oal)\in\Z$ or (ii) $\si(\al)=-\al$, in this case
$c_\al=4$ and
$(\mu,\oal)/(\oal,\oal)\in\Z$ if and only if $\langle\mu,\al^\vee\rangle\in2\Z$.

Being $\tPhi$ a root system and using the first statement of the theorem,
it follows that a special weight $\mu$ is spherical if and only if, for
$i=1,\ldots,\ell$, we have
$$
\left\{
\begin{array}{ll}
\langle\mu,\tal_i^\vee\rangle\in\Z & \textrm{ if } 2\tal_i\not\in\tPhi,\\
\langle\mu,(2\tal_i)^\vee\rangle\in\Z & \textrm{ if } 2\tal_i\in\tPhi,
\end{array}
\right.
$$
since $\tDe$ is a base for $\tPhi$.

If we assume that $\tPhi$ is reduced then the description of the generators
of $\Om_1$ follows from the above discussion. So let $\tPhi$ be non reduced.
We can suppose that $\tPhi$ is irreducible; it follows that $\tPhi$
is of type $\typeBC_\ell$.

First we consider the case of $X$ exceptional and let $\al\doteq\al_i$ be an exceptional root. We have
$$
\be\doteq-s_{\al}\si(\al)=-\si(\al)+\langle\si(\al),\al^\vee\rangle\al=-\si(\al)+\al
$$
by Lemma \ref{lexceptionalval}. Observe that $\tbe=2\tal$. So $\tal$ is the unique simple root in $\typeBC_\ell$ such that $2\tal\in\tPhi$. Let $\mu\doteq\om_i+\om_{\si(i)}$, we have
\begin{align*}
\langle\mu,(2\tal)^\vee\rangle
 & = \frac{(\mu,\overline{2\al})}{(\overline{2\al},\overline{2\al})}
   = \frac{1}{2}\frac{(\mu,\oal)}{(\oal,\oal)} =\\
 & = \frac{\langle\mu,\al_i^\vee\rangle}{c_\al}
   = \langle\mu,\al_i^\vee\rangle = 1,
\end{align*}
using again Lemma \ref{lexceptionalval}. Hence the result about $\Om_1$
holds also in this case.

If $X$ is non exceptional (and $\tPhi$ is non reduced) then the involution is
described, up to isomorphism, by one of the following two Satake diagrams
$$
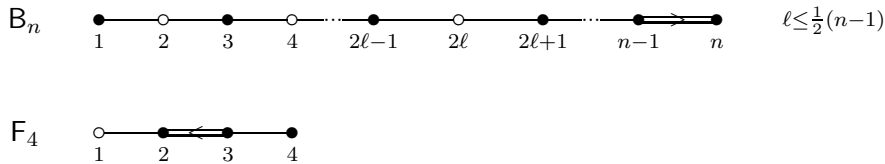
\hfil
\xymatrix@=4ex{
 *{\underset{\mathstrut}{\typeB_n}}
&*{\underset{\mathstrut 1}{\bullet}}
&*{\underset{\mathstrut 2}{\circ}}   \ar @{-} @<-1.2ex>  []-<0.4ex,0ex>;[l]+<0.4ex,0ex>
&*{\underset{\mathstrut 3}{\bullet}} \ar @{-} @<-1.2ex>  []-<0.4ex,0ex>;[l]+<0.4ex,0ex>
&*{\underset{\mathstrut 4}{\circ}} \ar @{-} @<-1.2ex>  []-<0.4ex,0ex>;[l]+<0.4ex,0ex>
&*{\underset{2\ell \mathstrut - 1}{\bullet}}
\ar @{-} @<-1.2ex> []-<0.4ex,0ex>;[l]+<0.4ex,0ex>|{\ldots}
&*{\underset{2 \mathstrut \ell}{\circ}}      \ar @{-} @<-1.2ex>  []-<0.4ex,0ex>;[l]+<0.4ex,0ex>
&*{\underset{2\ell \mathstrut + 1}{\bullet}} \ar @{-} @<-1.2ex>  []-<0.4ex,0ex>;[l]+<0.4ex,0ex>
&*{\underset{n\mathstrut - 1}{\bullet}}
\ar @{-} @<-1.2ex> []-<0.4ex,0ex>;[l]+<0.4ex,0ex>|{\ldots}
&*{\underset{\mathstrut n}{\bullet}}
\ar @{=} @<-1.2ex>  []+0;[l]+0 \ar @{}  @<-1.2ex>  [l] |{>}
& {\underset{\mathstrut}{\scriptstyle{\ell\leq\frac{1}{2}(n-1)}}}\\
 *{\underset{\mathstrut}{\typeF_4}}
&*{\underset{\mathstrut 1}{\circ}}
&*{\underset{\mathstrut 2}{\bullet}} \ar @{-} @<-1.2ex>  []+0;[l]+<0.4ex,0ex>
&*{\underset{\mathstrut 3}{\bullet}} \ar @{=} @<-1.2ex>  []+0;[l]+0 
                                     \ar @{}  @<-1.2ex>  [l] |{<}
&*{\underset{\mathstrut 4}{\bullet}} \ar @{-} @<-1.2ex>  []+0;[l]+0 \\
}
$$
In the first case we have
$$
\si(\al_i)=
\left\{
\begin{array}{ll}
-\al_i-(\al_{i-1}+\al_{i+1}) & \textrm{if } 2|i \textrm{ and } i<2\ell,\\
-\al_{2\ell}-(\al_{2\ell-1}+2\al_{2\ell+1}+\cdots+2\al_{n-1}+\al_n) &
\textrm{if } i=2\ell,\\
\al_i & \textrm{otherwise.}
\end{array}
\right.
$$
Since $\langle\si(\al_{2\ell}),\al_{2\ell}^\vee\rangle=1$ we have: (i) $s_{\al_{2\ell}}\si(\al_{2\ell})=\si(\al_{2\ell})-\al_{2\ell}$,
(ii) $(\tal_{2\ell},\tal_{2\ell})=(\al_{2\ell},\al_{2\ell})$.
Hence $\al_{2\ell}$ is the unique simple root such that
$\tal_{2\ell},2\tal_{2\ell}\in\tPhi$ and
\begin{align*}
\langle\om_{2\ell},(2\tal_{2\ell})^\vee\rangle &
   = \frac{(\om_{2\ell},\tal_{2\ell})}{(\tal_{2\ell},\tal_{2\ell})}
   = \frac{(\om_{2\ell},\al_{2\ell}-\si(\al_{2\ell}))}
{(\al_{2\ell},\al_{2\ell})} =\\
 & = \langle\om_{2\ell},\al_{2\ell}^\vee\rangle = 1.
\end{align*}
So the claimed description of $\Om_1$ is proved.

In the second case we have
$$
\si(\al_i)=
\left\{
\begin{array}{ll}
-\al_1-(3\al_2+2\al_3+\al_4) & \textrm{if } i=1,\\
\al_i & \textrm{otherwise.}
\end{array}
\right.
$$
Since $\langle\si(\al_1),\al_1^\vee\rangle=1$ we conclude
as in the previous case.
\end{proof}
\begin{corollary}\label{cdominantchamber}
\emph{1)}
The set $\Om_1^+$ of integral weights of $\tPhi$ that are dominant with
respect to $\tDe$ is $\La^+\cap\Om_1$;\\
\emph{2)} there exists a $\Z$ basis $\theta_1,\ldots,\theta_r$ for $\Pic(X)$ that is a $\N$ basis
for the cone $\Pic^+(X)$; moreover any $\theta_i$ is of the following three kind: $w_j$ or $2w_j$ or $w_j+w_{\si(j)}$ for some $1\leq j\leq n$.
\end{corollary}
\begin{proof}
It follows from the formula for the weights $\tom_i$ given in Theorem
\ref{tspherical} above and from Proposition \ref{ppic}.
\end{proof}

We introduce a new order on the set of weights.
Given two integral weights $\mu,\la\in\La$ we write $\mu\leq_\si\la$ if $\mu=\la-\tal$ for some $\tal\in\tDe_\N$, where $\tDe_\N$ is the positive cone over $\tDe$.

Let $\la\in\Pic^+(X)$. Using the order $\leq_\si$ above and Proposition
\ref{pfiltration} we can decompose $H^0(X,\lb_\la)$ as
$$
H^0(X,\lb_\la)=\oplus V_\mu^*,
$$
where the sum runs over the dominant weights $\mu$ such that $\mu\leq_\si\la$.
Now recall the PRV conjecture proved independently in \cite{KU} and \cite{M}.
\begin{proposition}\label{pprv}
Let $\la$ and $\mu$ be two dominant weights and let $\tau,\ep$ be two elements in the Weyl group. If $\nu=\tau(\la)+\ep(\mu)$ is dominant then the module $V_\nu$ appears in the decomposition into irreducible modules of the tensor product $V_{\la}\otimes V_{\mu}$.
\end{proposition}
A usefull consequence is the following
\begin{proposition}\label{pkprv}\emph{(Lemma 3.2 in \cite{K})}
Let $\nu,\la,\mu$ be dominant weights such that $\nu\leq\la+\mu$. Then there exist dominant weights $\la',\mu'$ such that (i) $\la'\leq\la$, $\mu'\leq\mu$ and (ii) $V_\nu$ appears in the decomposition of $V_{\la'}\otimes V_{\mu'}$.
\end{proposition}
We want to prove the following generalization for symmetric varieties
\begin{theorem}\label{tsprv}
If $\nu,\la,\mu\in\Om_1^+$ and $\nu\leq_\si\la+\mu$, then there exist two
weights $\la',\mu'\in\Om_1^+$ such that
(i) $\la'\leq_\si\la$, $\mu'\leq_\si\mu$ and
(ii) $V_\nu$ appears in the decomposition of $V_{\la'}\otimes V_{\mu'}$.
\end{theorem}
\begin{proof}
For $\eta\in\lit^*$ we denote by $[\eta]_W$ the unique element of the dominant Weyl chamber in the orbit $W\cdot\eta$. For $\eta\in\lit_1^*$ we define analogously $[\eta]_{\tW}$ using the root system $(\tPhi,E_1)$. We observe that if $\eta\in\lit_1^*$ then $[\eta]_W=[\eta]_{\tW}$ since by Proposition \ref{prestricted} we have $\tW\cdot\eta\subset W\cdot\eta$ and since the dominant Weyl chamber of $\tPhi$ is contained in the dominant Weyl chamber of $\Phi$ by Corollary \ref{cdominantchamber}.

Assume that $\tPhi$ is reduced and let $K$ be a simply connected group with root system $\tPhi$. Notice that the set $\Om_1^+$ is the set of dominant integral weights of $K$ and $\leq_\si$ is the dominant order for $K$. For $\eta\in\Om_1^+$ let $\Om_\eta$ be the set of weights of the irreducible $K$ module $\tV_\eta$ of highest weight $\eta$.

We have $\nu\in\Om_{\la+\mu}$ since $\nu\leq_\si\la+\mu$ and $\nu,\la,\mu\in\Om_1^+$. Consider now the $K$ equivariant projection
$$
\tV_\la\otimes\tV_\mu\twoheadrightarrow\tV_{\la+\mu}.
$$
We have that there exist weights $\ola\in\Om_\la$, $\omu\in\Om_\mu$ such that $\nu=\ola+\omu$. In particular $\la-\ola\in\tDe_\N$, $\mu-\omu\in\tDe_\N$. Let $\la'\doteq[\ola]_W$, $\mu'\doteq[\omu]_W$ and notice that $\la',\mu'\in\Om_1^+$ by the remark at the beginning of the proof. Moreover $\la'-\ola,\mu'-\omu\in\tDe_\N$ since $\ola,\omu$ are integral weights. So $\la-\la',\mu-\mu'\in\tDe_\Z$. Further $\la'\leq\la$, $\mu'\leq\mu$ since $\la'\in\Om_\la$, $\mu'\in\Om_\mu$. We conclude $\la'\leq_\si\la$, $\mu'\leq_\si\mu$. This shows also $\la',\mu'\in\Om_1^+$. Finally $V_\nu$ appears in $V_{\la'}\otimes V_{\mu'}$ using the PRV (Proposition \ref{pprv}).

Assume now that $\tPhi$ is non reduced and $\tPhi$ is irreducible (without loss of generality); so $\tPhi$ is of type $\typeBC_\ell$. The proof given above still holds with the following remarks: (i) choose $K$ of type $\typeB_\ell\subset\typeBC_\ell$ and notice that $\leq_\si$ is the dominant order for $\typeB_\ell$; (ii) $\Om_1^+$ is strictly contained in the set of dominant integral weights of $K$ and (iii) if $\eta\in\Om_1^+$ and $\zeta\leq_\si\eta$ then $\eta\in\Om_1^+$.
\end{proof}
As an application to symmetric varieties we have
\begin{corollary}\label{csurcomb}
Let $X$ be a non exceptional complete symmetric variety and let $\la,\mu\in\Pic^+(X)$. Consider the multiplication map
$$
H^0(X,\lb_\la)\otimes H^0(X,\lb_\mu)\rightarrow H^0(X,\lb_{\la+\mu}).
$$
Then if $V_\nu^*$ appears in the right hand side as a direct summand then it appears in the left hand side too.
\end{corollary}
\begin{proof}
Follows from Theorem \ref{tsprv} above and Proposition \ref{ppic}.
\end{proof}

\section{The Standard Monomial Theory}\label{ssmt}

We begin this section with a short review of Littelmann path basis theory, for details see \cite{L1},\cite{L2} and \cite{LLM}. Denote by $\Pi$ the set of piecewise linear paths $\pi:[0,1]_\Q\rightarrow\La\otimes\Q$ starting in $0$. Let $\Pi^+$ be the set of path $\pi\in\Pi$ such that $\img\pi$ is contained in the dominant Weyl chamber. For $\al\in\De$ let $e_\al,f_\al$ be the root operators associated with $\al$. We can associate a coloured directed graph $\G(\B)$ to any subset $\B$ of $\Pi$ by joining two paths $\pi_1,\pi_2\in\B$ with an arrow $\xymatrix@1{\pi_1\ar[r]^\al & \pi_2}$ if $f_\al(\pi_1)=\pi_2$.

Now let $\la\in\La^+$ and choose $\pi\in\Pi^+$ such that $\pi(1)=\la$, then the \emph{path model} associated with $\pi$ is the set $\B_\pi$ of paths obtained from $\pi$ by applying the root operators, i.e., $\B_\pi\cup\{0\}$ is the smallest subset of $\Pi\cup\{0\}$ that contains $\pi$ and is closed under the root operators $e_\al, f_\al$. The path model describes the character of the $\cG$ module $V_\la$, indeed we have $\ch V_\la=\sum_{\eta\in\B_\pi}e^{\eta(1)}$.

In particular if we start with $\pi_\la:t\mapsto t\la$ then $\B_\la\doteq\B_{\pi_\la}$ is the set of Lakshmibai-Seshadri paths (LS paths) of shape $\la$. This path model has a simple combinatorial description in terms of poset with bonds (see \cite{CHI1}, \cite{CHI2}). Let $W_\la$ be the stabilizer of $\la$ in the Weyl group $W$ and consider the set of minimal representatives $W^\la$. Given two adjacent elements $\tau_1<s_\al\tau_2$ in $W^\la$, where $\al\in\De$, one can define a positive integer value function $f_\la$ as $f_\la(\tau_1,\tau_2)=\langle\tau_1(\la),\al^\vee\rangle$. Further given a complete chain $\tau_1<\cdots<\tau_u$ in $W^\la$ we have that $\gcd\{f_\la(\tau_1,\tau_2),\ldots,f_\la(\tau_{u-1},\tau_u)\}$ depends only on the pair $\tau_1,\tau_u$; so one can extend $f_\la$ to comparable pairs defining $f_\la(\tau_1,\tau_u)=\gcd\{f_\la(\tau_1,\tau_2),\ldots,f_\la(\tau_{u-1},\tau_u)\}$. The data $(W^\la,\leq,f_\la)$ is called a \emph{poset with bonds}.

Then the set $\B_\la$ of LS paths of shape $\la$ is in bijection with the set of pairs $(\tau_1<\cdots<\tau_u;0=a_0<a_1<\cdots<a_{u-1}<a_u=1)$ such that $a_if_\la(\tau_i,\tau_{i+1})\in\N$ for $i=1,\ldots,u-1$. Now let $N_\la$ be the least common multiple of the image $f_\la(W^\la)$ and consider the set $\W$ of words in the alphabet $W$. We define the \emph{word} $w(\pi)\in\W$ of an LS path $\pi=(\tau_1<\cdots<\tau_u;0=a_0<a_1<\cdots<a_u=1)$ as
$w(\pi)=\tau_1^{N_\la(a_1-a_0)}\cdots\tau_u^{N_\la(a_u-a_{u-1})}$.

Now let $\la_1,\ldots,\la_r$ be dominant weights and call \emph{formal monomial} any string $\pi_1\cdots\pi_r$ with $\pi_i\in\B_{\la_i}$ for $i=1,\ldots,r$. In order to introduce a variant of the lexicographic order, we extend the word to monomials defining $w(\pi_1\cdots\pi_u)=w(\pi_1)\cdots w(\pi_u)$ using juxtaposition of words in $\W$.
Given two monomials $\pi_1\cdots\pi_r$, $\eta_1\cdots\eta_r$ we define
$\pi_1\cdots\pi_r\leq\eta_1\cdots\eta_r$ if
$w(\pi_1\cdots\pi_r)\lex\tau_1(w(\eta_1))\cdots\tau_r(w(\eta_r))$ for all permutations
$(\tau_1,\ldots,\tau_r)$ in the subgroup $\perm_{N_{\la_1}}\times\cdots\times\perm_{N_{\la_r}}$ of the symmetric group $\perm_{N_{\la_1}+\cdots+N_{\la_r}}$ on $N_{\la_1}+\cdots+N_{\la_r}$ symbols acting in the natural way on words of length $N_{\la_1}+\cdots+N_{\la_r}$.

Consider the set $\B_{\la_1}*\cdots*\B_{\la_r}$ of concatenations of LS paths of shapes $\la_1,\ldots,\la_r$ and notice that it is stable under the root operators. Define $\G(\la_1,\ldots,\la_r)$ as the connected component of the graph $\G(\B_{\la_1}*\cdots*\B_{\la_r})$ containing $\pi_{\la_1}*\cdots*\pi_{\la_r}$; let $\la\doteq\la_1+\cdots+\la_r$ and recall that the map $\pi_{\la_1}*\cdots*\pi_{\la_r}\mapsto\pi_\la$ extends to an isomorphism of graphs $\G(\la_1,\ldots,\la_r)\rightarrow\G(\B_\la)$. Now given a formal LS path monomial $\pi_1\cdots\pi_r$ we call it \emph{standard} if $\pi_1*\cdots*\pi_r\in\G(\la_1,\ldots\la_r)$. It is then clear that the number of formal standard monomial $\pi_1\cdots\pi_r$ is given by $|\G(\la_1,\ldots,\la_r)|=|\B_\la|=\dim V_\la$.

This machinery has another key feature for our purpose. 
Let $\la_1,\ldots,\la_r$ be dominant weights such that the 
lattice $\Psi\doteq\langle\la_1,\ldots,\la_r\rangle_\Z$ has rank $r$ and 
let $Q$ 
be the parabolic subgroup $P_{\lambda_1 + \dots + \lambda_r}$ of $G$. 
Consider the ring
$A_\Psi(G/Q)\doteq\oplus_{\la\in\Psi}H^0(G/Q,\lb_\la)$.
Let $\la\in\Psi\cap\La^+$ and recall that Littelmann associates a 
section $p_\pi\in H^0(G/Q,\lb_\la)$ to an LS path $\pi\in\B_\la$. These 
sections have very 
remarkable properties and they give a standard monomial theory for 
$A_\Psi(G/Q)$. Let $\pi_1,\ldots,\pi_u$ be LS paths with 
$\pi_i\in\B_{\la_{h_i}}$ and $h_1\leq\cdots\leq h_u$, and call a monomial 
$p_{\pi_1}\cdots p_{\pi_u}$ standard if $\pi_1\cdots\pi_u$ is a standard 
formal monomial, further we define the shape of the monomial 
$p_{\pi_1}\cdots p_{\pi_u}$ as $\la_{h_1}+\cdots+\la_{h_u}$. 
Then one has
\begin{proposition}\label{pflagsmt}
\emph{1)} The standard monomials of shape $\la$ forms a $\field$ basis for $H^0(G/Q,\lb_\la)$;\\
\emph{2)} if $p_{\pi_1}\cdots p_{\pi_u}=\sum_h a_h p_{\eta_{h,1}}\cdots p_{\eta_{h,u}}$
express the non standard monomial $p_{\pi_1}\cdots p_{\pi_u}$ in terms of standard monomials with $a_h\neq0$ then $\eta_{h,1}\cdots\eta_{h,u}\leq\pi_1\cdots\pi_u$ for all $h$.
\end{proposition}
We will refer to the expressions in 2) of Proposition \ref{pflagsmt} above as the \emph{Littelmann relations}.

We are ready to develop our standard monomial theory. First a definition, for a variety $Z$ let $A(Z)$ denote the ring $\oplus_{\lb\in\Pic(Z)}H^0(Z,\lb)$ that we call the \emph{ring of sections} of $Z$. Now let $X$ be a complete symmetric variety and let $\theta_1,\ldots,\theta_r$ be the $\Z$ basis for $\Pic(X)$ as in Corollary \ref{cdominantchamber}. Let $\B_i\doteq\B_{\theta_i}$ be the LS path basis of shape $\theta_i$ described above for $i=1,\ldots,r$. For an LS path $\pi\in\B_i$ it is
possible, according to Proposition \ref{pfiltration}, to choose $x_\pi\in H^0(X,\lb_{\theta_i})$ such that $x_{\pi|G/P}=p_\pi$.
\begin{lemma}\label{lgenerators}
The sections $x_\pi$ for $\pi\in\B_1\sqcup\cdots\sqcup\B_r$ and the sections
$s_i\in H^0(X,\lb_{\tal_i})$ for $i=1,\ldots,\ell$ generate $A(X)$.
\end{lemma}
\begin{proof}
Let $A'(X)$ be the $\field$ subalgebra generated by the sections in the
statement. By Proposition \ref{pfiltration} it is enough to show that for
each dominant $\la\in\Pic(X)$ the submodule $V_\la^*\subset H^0(X,\lb_\la)$
is contained in $A'(X)$.

For a weight $\la=a_1\theta_1+\cdots+a_r\theta_r$ let
$\he\la\doteq\sum_i a_i$. We use induction on $\he\la$. If $\he\la=0$ then $\la=0$ and $H^0(X,\lb_\la)\simeq\field\cdot1$. So suppose $\he\la>0$.
Hence there exists a dominant weight $\la'$ and $1\leq i\leq r$ such that $\la=\la'+\theta_i$. Consider the following commutative diagram
\begin{figure}[ht!]
\hfil
\xymatrix{
H^0(X,\lb_{\la'})\otimes H^0(X,\lb_{\theta_i})\ar[r] & H^0(X,\lb_\la)\\
V_{\la'}^*\otimes V_{\theta_i}^*\ar@{^{(}->}[]+<0ex,2.5ex>;[u]\ar@{>>}[r]
& V_\la^*\ar@{^{(}->}[]+<0ex,2.5ex>;[u]
}
\hfil
\end{figure}
where the horizontal maps are induced by multiplication. Notice that the
lower horizontal map is surjective by $\cG$ equivariance.
This finish the proof since $\he\la'<\he\la$.
\end{proof}
A generic monomial in the generators can be written in the form
$$
x=s_1^{n_1}\cdots s_\ell^{n_\ell}x_{\pi_1}\cdots x_{\pi_u}
$$
with $\pi_i\in\B_{h_i}$ and $h_1\leq\cdots\leq h_u$. For such a monomial
$x$ we define the \emph{order of vanishing} as $v(x)\doteq(n_1,\ldots,n_\ell)$,
the \emph{shape} as
$\la(x)\doteq\sum_{i=1}^\ell n_i\tal_i+\sum_{j=1}^h\theta_{h_j}$ and the
\emph{flag shape} as $\mu(x)\doteq\la(\ox)$ where $\ox\doteq s^{-v(x)}x$.
Notice that $x\in F_{\la(x)}(v(x))\subset H^0(X,\lb_{\la(x)})$ and that
$\mu(x)\in\La^+$.

We define the set $\M$ of \emph{standard monomials} for $X$ as the set of
monomials $x=s_1^{n_1}\cdots s_\ell^{n_\ell}x_{\pi_1}\cdots x_{\pi_u}$ as above
such that $\pi_1\cdots\pi_u$ is a formal LS standard monomial. Further we denote by $\M_\la$ the set of standard monomials $x$ such that $\la(x)=\la$.
Given two monomials
$x=s_1^{n_1}\cdots s_\ell^{n_\ell}x_{\pi_1}\cdots x_{\pi_u}$ and
$y=s_1^{m_1}\cdots s_\ell^{m_\ell}x_{\eta_1}\cdots x_{\eta_v}$ with the same shape we write
$x\leq y$ if $v(x)<v(y)$ or $v(x)=v(y)$ and $\pi_1\cdots\pi_u\leq\eta_1\cdots\eta_v$.

Finally let $\A(X)$ be the polynomial ring with indeterminates
$s_1,\ldots,s_\ell$ and $x_\pi$ with $\pi\in\B_1\sqcup\cdots\sqcup\B_r$.
Clearly $A(X)$ is isomorphic to a quotient of $\A(X)/I$ for some ideal $I\subset\A(X)$.

The main result of our standard monomial theory is the following
\begin{theorem}\label{tsmt}
\emph{1)} The set $\M_\la$ is a $\field$ basis for $H^0(X,\lb_\la)$;\\
\emph{2)} given monomials $x_1,\ldots,x_t\in\M$, let
$x_1\cdots x_t=\sum a_z z$ with $z\in\M$ be the relation guaranteed by
\emph{1)}. Then for any standard monomial $z$ such that $a_z\neq0$ we have
$x_1\cdots x_t\leq z$. Moreover $\ox_1\cdots\ox_t=\sum a_z\oz$ with
$v(z)=v(x)+v(y)$ is the Littelmann relation for a multicone over $G/P$;\\
\emph{3)} the ideal $I$ is generated by the relation in \emph{2)} for $t=2$
and $x=x_{\pi_1}$, $x_2=x_{\pi_2}$ with
$\pi_1,\pi_2\in\B_1\sqcup\cdots\sqcup\B_r$ and $x_{\pi_1}x_{\pi_2}$ not
standard.
\end{theorem}
\begin{proof}
We prove the three statements together.

First notice that a section in $A(X)$ vanishing on $G/P$ is in the ideal
generated by $s_1,\ldots,s_\ell$ since the divisor $S_1,\ldots,S_\ell$ are smooth and normal crossing.

If $\pi_1,\pi_2\in\B_1\sqcup\cdots\sqcup\B_r$ are two LS paths such that
$x_{\pi_1}x_{\pi_2}$ is not standard, consider the Littelmann relation
$$
p_{\pi_1}p_{\pi_2}=\sum_h a_h p_{\eta_{h,1}}p_{\eta_{h,2}}
$$
on $G/P$. Then $x_{\pi_1}x_{\pi_2}-\sum_h a_h x_{\eta_{h,1}}x_{\eta_{h,2}}$
vanishes on $G/P$, hence
$$
x_{\pi_1}x_{\pi_2}=\sum_h a_h x_{\eta_{h,1}}x_{\eta_{h,2}}+
\sum_{v(z)>0}a_z z
$$
where in the second sum the $z$'s are (not necessarily standard) monomials.

Consider the $\La$ homogeneous element
$$
f_{\pi_1,\pi_2}\doteq x_{\pi_1}x_{\pi_2}-
\sum_h a_h x_{\eta_{h,1}}x_{\eta_{h,2}}-\sum_{v(z)>0}a_z z
$$
in $\A(X)$ and let $J$ be the ideal generated by the various $f_{\pi_1,\pi_2}$
with $x_{\pi_1}x_{\pi_2}$ not standard as above. We want to show that $\M$
generates $\A(X)/J$ as a vector space.

Let $x\doteq s_1^{n_1}\cdots s_\ell^{n_\ell}x_{\pi_1}\cdots x_{\pi_u}$ be a not
standard monomial (so $u\geq2$). We proceed by induction on the flag shape
of $x$ with respect to the order $\leq_\si$.

Observe that $\A(X)/(\langle s_1,\ldots s_\ell\rangle+J)$ is isomorphic to the
coordinate ring of a multicone over $G/P$ since the relations on sections over $G/P$ are generated by the relations of degree 2 (Proposition 2 in \cite{KR}). Hence $x_{\pi_1}\cdots x_{\pi_u}+\langle s_1,\ldots,s_\ell\rangle$ is a sum
of standard monomials in the $x_\pi$. So in $\A(X)/J$ we have
$x_{\pi_1}\cdots x_{\pi_u}=\sum_{z\in\M}a_z z+s_1 y_1+\cdots+s_\ell y_\ell$,
where $y_1,\ldots,y_\ell$ are sums of monomials with flag shape $<_\si$ of
the flag shape of $x$.

Now consider the $\La$ homogeneous projection $\phi:\A(X)/J\rightarrow A(X)$.
We want to show that $\phi$ is an isomorphism. It is enough to prove
that $\dim(\A(X)/J)_\la\leq\dim(A(X))_\la$ for each $\la$, since $\phi$ is
clearly surjective. We have $\dim(\A(X)/J)_\la\leq|\M_\la|$. On the other hand
$(A(X))_\la=H^0(X,\lb_\la)=\oplus V_\mu^*$ with $\mu\leq_\si\la$, $\mu$
dominant. So
$\dim(A(X))_\la=\sum_{\mu\leq_\si\la,\ \mu\in\La^+}\dim V_\mu^*=
\sum_{\mu\leq_\si\la,\ \mu\in\La^+}|\B_{\pi_\mu}|$.
If $\mu=a_1\theta_1+\cdots+a_r\theta_r$ then
$$
|B_\mu|=|\G(
\underbrace{\theta_1,\ldots,\theta_1}_{a_1},\ldots,\underbrace{\theta_r,\ldots,\theta_r}_{a_r}
)|
$$
and we conclude $\sum_{\mu\leq_\si\la,\ \mu\in\La^+}|\B_\mu|=|\M_\la|$.
\end{proof}

Now we show that this standard monomial theory is compatible with the $G$ orbit closures in $X$. So let $I\subset\{1,\ldots,\ell\}$ and let $X_I\doteq\cap_{i\in I}S_i$ be the corresponding $G$ orbit closure. Define a monomial $x=s_1^{n_1}\cdots s_\ell^{n_\ell}x_{\pi_1}\cdots x_{\pi_u}$ to be \emph{standard on $X_I$} if it is standard and $n_i=0$ for all $i\in I$. Given a weight $\la$, denote by $\M_{\la,I}$ the set of monomials standard on $X_I$ with shape $\la$. Then we have
\begin{corollary}\label{cgvarietybasis}
The set $\M_{\la,I}$ is a $\field$ basis for $H^0(X_I,\lb_{\la|X_I})$.
\end{corollary}
\begin{proof}
Let $J$ be the complement of $I$ in $\{1,\ldots,\ell\}$. Adapting the proof of Theorem~8.3 in \cite{CP} we have that $H^0(X,\lb_\la)\rightarrow H^0(X_I,\lb_{\la|X_I})$ is a surjective map and that $H^0(X_I,\lb_{\la|X_I})=\oplus_\mu V_\mu^*$ where the sum runs over all dominant weights $\mu$ of the form
$$
\mu=\la-\sum_{j\in J}a_j\tal_j
$$
with $a_j\geq 0$. So the set of standard monomials $\M_\la$ is a generating set for $H^0(X_I,\lb_{\la|X_I})$. Hence $\M_{\la,I}$ is a generating set since any monomial in $\M_\la\setminus\M_{\la,I}$ contains some $s_i$ with $i\in I$ and vanishes on $X_I$. Moreover comparing the dimensions we have that $\M_{\la,I}$ is a basis.
\end{proof}
Clearly all the statements of Theorem \ref{tsmt} carry on to $X_I$ giving a standard monomial theory for the ring $A_{\Pic(X)}(X_I)=\oplus_{\la\in\Pic(X)}H^0(X_I,\lb_{\la|X_I})$.

Now we want to give some straightforward applications of the discussion above. 
The form of the relations in Theorem \ref{tsmt} allows to degenerate
the ring $A(X)$ to the coordinate ring of the product of a multicone
over the flag varieties $G/P$ and the affine space $\field^\ell$. Indeed, 
let $K$ be the ideal of $A(X)$ generated by $s_1,\ldots,s_\ell$ and consider 
the Rees algebra
$$
\Rees\doteq\cdots\oplus A(X)t^2\oplus A(X)t\oplus A(X)\oplus Kt^{-1}\oplus
K^{2}t^{-2}\oplus\cdots\subset A(X)\otimes\field[t,t^{-1}].
$$
Let $\tA(G/P)\doteq A_{\Pic(X)}(G/P)$. Then we have
\begin{theorem}\label{tdegeneration}
$\Rees$ is a flat $\field[t]$ algebra. The general fiber $\Rees/(t-a)$, with $a\in\field\setminus\{0\}$, is isomorphic to $A(X)$ and the special fiber $\Rees/(t)$ is isomorphic to $\tA(G/P)\otimes\field[s_1,\ldots,s_\ell]$.
\end{theorem}
\begin{proof}
This is a standard result about Rees algebra taking into account the relations of Theorem \ref{tsmt}.
\end{proof}
Now we use this degeneration result to prove that $A(X)$ has rational singularities.
As we said in the introduction, this is well known and we include it here since we haven't found it in the literature.
In  the proof below we will use (i) that $A(G/P)$ has rational singularities (Theorem 2 in \cite{KR}) and (ii) that the property of having rational singularities is stable under flat deformation (see \cite{E}) and under the quotient by a reductive group (see \cite{B}).
\begin{theorem}\label{trational}
\emph{1)} The ring $A(X)$ has rational singularities;\\
\emph{2)} for all $\lb\in\Pic(X)$ the ring $A_\lb\doteq\oplus_{n\in\N}H^0(X,\lb^{\otimes n})$ has rational singularities.
\end{theorem}
\begin{proof}
By what we said above and Theorem \ref{tdegeneration} in order to prove 1) it is enough to show that $\tA(G/P)$ is the fixed point algebra of $A(G/P)$ under the action of a reductive group. This will be done in two steps. 
Let $\Psi$ be the lattice 
$\langle \theta_1,\dots,\theta_r\rangle_{\Q} \cap \Lambda$ and observe that 
$$
\Lambda \supset 
\Pic(G/P)= \langle \omega_{\alpha} \colon \alpha \in \Delta_1 \rangle_{\Z} \supset
\Psi \supset \Pic(X)
$$
Let $\tT$ be the maximal torus of $\cG$ over $T$ and recall that
$\Lambda = \Hom (\tT , \field^*)$. We define $\tS\doteq\bigcap_{i=1}^{r} \ker 
\theta_i \subset \tT$ and an action of $\tS$ on $A(G/P)$ by
$$
t\cdot f\doteq\mu(t)f \textrm{ for all } t\in\tS \textrm{ and } 
f\in H^0(G/P,\lb_\mu).
$$
We have that $(A(G/P))^{\tS}= A_\Psi(G/P)$.
Now we observe that $\Psi/ \Pic (X)$ is a finite abelian group.
Hence $\Gamma = \Hom(\Psi/\Pic(X) , \field^*)$ is also finite and 
we can define an action of $\Gamma$ on $A_{\Psi}(G/P)$ by
$$
\ga\cdot f\doteq\ga(\mu+\Psi)f \textrm{ for all } \ga\in\Ga \textrm{ and } f\in H^0(G/P,\lb_\mu).
$$
Clearly $(A_{\Psi}(G/P))^\Ga = \tA(G/P)$.

For the proof of 2) notice that if $\lb$ is not trivial and there exists $m>0$ such that $H^0(X,\lb^{\otimes m})\neq0$ then $H^0(X,\lb^{\otimes n})=0$ for all $n<0$. So we can assume $A_\lb(X)=\oplus_{n\in\Z}H^0(X,\lb^{\otimes n})$
and $\lb = \lb_{\lambda}$ for some $\lambda \in \Lambda$. 
Then the proof goes on as in 1) using the action of a
torus and a finite group to pass from the algebra $A(X)$, 
corresponding to $\Pic (X)$, to $A_\lb(X)$, corresponding to $\Z\la$.
\end{proof}

As a final remark we consider the characteristic $p$ case.
In this paper we have treated exclusively the characteristic $0$ case. In
particular in the proof of the last Theorem we used the result of Boutot 
on quotients of rational singularities and the result of Elkike on
deformation of rational singularities which in general hold only in this
case and we do not know if Theorem \ref{trational} holds also in the 
finite characteristic case. 

However it is possible to prove Theorem \ref{tsmt} in general with few changes to our proof.
We give now an outline of the modification needed to pass
from the characteristic zero case to the general case.

The characterizations of $\Omega_1$ and $\Pic (X)$ given in Theorem 
\ref{tspherical} and Corollary \ref{cdominantchamber} hold in the same
way and with the same proofs in the finite characteristic case. The main
changes are in the proof of Lemma \ref{lgenerators}. One minor change is
that we cannot use the decomposition of $H^0(X,\lb_\la) =
\bigoplus_{\mu\leq_\si\la,\mu\in\La^+} V_{\mu}^*$, but we have instead to
use the filtration $F_\la$ of Proposition \ref{pfiltration}. 
A more important change is that $V_\la^*$ is no more irreducible but we can instead use the surjectivity of the multiplication of sections
$H^0(G/P,\lb_\la)\otimes H^0(G/P,\lb_\mu)\longrightarrow H^0(G/P,\lb_{\la+\mu})$ proved
by Ramanan and Ramanathan (see \cite{RR}). Then the remaining part of the proof of Theorem \ref{tsmt} goes on with only minor changes.

\vskip 0.5cm
\scriptsize{
Rocco Chiriv\`\i,\\
Universit\`a di Pisa,\\
Dipartimento di Matematica ``Leonida Tonelli''\\
via Buonarroti n.~2, 56127 Pisa, Italy,\\
e-mail: {\tt chirivi@dm.unipi.it}\\
\indent and\\
Andrea Maffei,\\
Universit\`a di Roma ``La Sapienza",\\
Dipartimento di Matematica ``Guido Castelnuovo",\\
Piazzale Aldo Moro n.~5, 00185 Rome, Italy,\\
e-mail: {\tt amaffei@mat.uniroma1.it}
}
\end{document}